\documentclass[12pt]{article}
\usepackage{amsfonts}
\usepackage{latexsym,amsmath,amssymb}
\usepackage{graphics,epstopdf}
\usepackage[pdftex]{graphicx}
\usepackage{subcaption}

\bibliography{bibfile}
\newtheorem{definition}{Definition}[section]
\newtheorem{theorem}[definition]{Theorem}

\newtheorem{lemma}[definition]{Lemma}
\numberwithin{equation}{section}

\newtheorem{thm}{Theorem}[section]

\newtheorem{rem}[thm]{Remark}
\numberwithin{equation}{section}

\begin{document}

	\bigskip
	
	\bigskip
	
	\begin{center}
		{\Large \textbf{ On King type modification of $(p,q)$-Lupa\c{s} Bernstein operators.
			}}
			
	\bigskip
			
	\textbf{ Asif Khan, Vinita Sharma }
	and \textbf{Faisal Khan}
			
	Department of\ Mathematics, Aligarh Muslim University,
	Aligarh 202002, India%
	\\[0pt]
			
	akhan.mm@amu.ac.in; vinita.sha23@gmail.com; faisalamu2011@gmail.com
			

	\bigskip
			
	\bigskip

	\textbf{Abstract}
	\end{center}
	\parindent=8mm {\footnotesize {In this paper, a King-type modification of $(p,q)$-Lupa\c{s} Bernstein operators are introduced. The rate of convergence of these operators are studied by  means of modulus of continuity and Lipschitz class functional.\\
	Further, it has been shown that the error estimation of these operators on some subintervals of $[0,1]$ are better than the $(p,q)$-Lupa\c{s} Bernstein operators.\\ 
	\bigskip
			

	{\footnotesize \emph{MSC: primary 65D17; secondary 41A10, 41A25, 41A36.}: \newline

 \section{Introduction}		
First, 	Let us recall certain notations of $(p,q)$-calculus.\\

For any $n\in \mathbb N$, the $(p,q)$-integers are defined as follows:
 	\begin{equation*}
 	\lbrack n]_{p,q}=p^{n-1}+p^{n-2}q+p^{n-3}q^2+...+pq^{n-2}+q^{n-1}\\
 	=\left\{
 	\begin{array}{lll}
 	\frac{p^{n}-q^{n}}{p-q},~~~~~~~~~~~~~~~~\mbox{when $~~p\neq q \neq 1$  } & \\
 	&  \\
 	n~p^{n-1},~~~~~~~~~~~~~~\mbox{ when $p=q\neq1$  } & \\
 	&  \\
 	
 	[n]_q ,~~~~~~~~~~~~~~~~~~~\mbox{when $p=1$  }& \\

 	n ,~~~~~~~~~~~~~~~~~~~~~\mbox{ when $p=q=1.$  }
 	\end{array}%
 	\right.
 	\end{equation*}
 Also the $(p,q)$-binomial coefficient is defined by\\

 $${n \brack k}_{p,q}=\frac{[n]_{p,q}!}{[k]_{p,q}!~[n-k]_{p,q}!}~~\text{for all}~~n, k\in \mathbb N~~\text{with}~~n\ge k.$$
		
When $ p=1$ and $q=1,$ it reduces to the ordinary integers and binomial cofficient respectively.\\

Recently, the applications of $(p,q)$-calculus emerged as a new area in the field of approximation theory. The $(p,q)$-calculus development has led to the discovery of various
			generalizations of Bernstein polynomials based on $(p,q)$-integers. The purpose of
			these generalizations is to provide appropriate and powerful tools to
			application areas such as computer-aided geometric
			design, numerical analysis,  and solutions of differential equations.\\
			
			Mursaleen \textit{et al}  \cite{mka1} first introduced $(p,q)$-calculus in approximation theory and constructed the {\ $(p,q)$-analogue of Bernstein operators defined as follows  for $0<q<p\leq 1$:\\

\begin{equation}\label{ee1}
B_{n,p,q}(f;x)=\frac1{p^{\frac{n(n-1)}2}}\sum\limits_{k=0}^{n}\left[
\begin{array}{c}
n \\
k%
\end{array}%
\right]_{p,q}p^{\frac{k(k-1)}{2}}x^{k}\prod \limits_{s=0}^{n-k-1}(p^{s}-q^{s}x)~~f\left( \frac{
[k]_{p,q}}{p^{k-n}[n]_{p,q}}\right) ,~~x\in \lbrack [0,1].
\end{equation}

			where

				\begin{align*}
				(1-x)^{n}_{p,q}&=\prod\limits_{s=0}^{n-1}(p^s-q^{s}x) =(1-x)(p-qx)(p^{2}-q^{2}x)...(p^{n-1}-q^{n-1}x)\\
				&=\sum\limits_{k=0}^{n} {(-1)}^{k}p^{\frac{(n-k)(n-k-1)}{2}} q^{\frac{k(k-1)}{2}}\left[
				\begin{array}{c}
				n \\
				k
				\end{array}%
				\right] _{p,q}x^{k}.
				\end{align*}

					Note when $p=1,$ $(p,q)$-Bernstein Operators given by $(\ref{ee1})$ turns out to be Phillips $q$-Bernstein Operators \cite{pl}.\\

The $q$-analogue of Bernstein operators \cite{sn} introduced by Lupa\c{s} \cite{lp} are as follows:\\

	\begin{equation}\label{e1.1}
	L_{n}(f;p;q;x)=\sum\limits_{k=0}^{n}f\biggl{(}\frac{k}{n}\biggl{)}b_{nk}(q;x),~~f\in C[0,1]~~ and~~  x\in [0,1].
	\end{equation}
	where\\
	
\begin{equation}\label{15e}
b_{nk}(q;x)=\frac{\left[\begin{array}{c}
	n \\
	k%
	\end{array}
	\right] _{q} q^{\frac{k(k-1)}{2}}~x^{k}~(1-x)^{n-k}}{\prod\limits_{j=0}^{n-1}\{(1-x)+q^{j} x\}}.
\end{equation}

Recently, Khalid et al. defined  $(p,q)$-analogue of Lupa\c{s} Bernstein operators \cite{khalidjcam}  as follows :\\

For any $p>0$ and $q>0,$ the linear operators $L^{n}_{p,q}:$ $C[0,1] \rightarrow C[0,1]$
\begin{equation}\label{e9a}
L^{n}_{p,q}(f;x)=L_{n}(f;p;q;x)= \sum\limits_{k=0}^{n}~~ \frac{f \bigg(\frac{p^{n-k}~[k]_{p,q}}{[n]_{p,q}}\bigg) ~\left[\begin{array}{c}
	n \\
	k%
	\end{array}%
	\right] _{p,q} p^{\frac{(n-k)(n-k-1)}{2}} q^{\frac{k(k-1)}{2}}~x^{k}~(1-x)^{n-k}}{\prod\limits_{j=1}^{n}\{p^{j-1}(1-x)+q^{j-1} x\}}.
\end{equation}

We recall the following lemma for the above operators.\\

\begin{lemma}\label{b1}\cite{khalidjcam}
	The following equalities are true\\
	
(1.)	$L_{n}(1;p;q;x)= 1$\\

(2.)	$L_{n}(t;p;q;x)= x$\\
	
(3.)	$L_{n}(t^{2};p;q;x)= x^2 +\frac{x(1-x)p^{n-1}}{[n]_{p,q}}~-~\frac{x^2(p-q)(1-x)}{p(1-x+qx}\big(1 - \frac{p^{n-1}}{[n]_{p,q}})$\\
\end{lemma}

Here, we know that the opeator $L_{n}(t^{2};p;q;x)$ do not preserve the test function $e_{2}$.\\

 In 2003, King \cite{king} introduced a non- trivial sequence of opeators preserving the functions $e_0$ and $e_2$ where $(e_i~=~x^i,i~=~0,1,2).$

He also proved that these operators have a better rate of convergence than the classical Bernstein polynomials whenever $0 \leq x \leq \frac{1}{3}$. In \cite{agritini}, Agratini and dogru introduced a King type modification of $q$-Szasz-Mirakjan type operators and they proved that their operators have a better rate of convergence than the classical
-Szasz-Mirakjan operators.\\

One can refer \cite{ kadak2,kadak5,kadak6,kanat,ogun,mishra2,mishra3,wn2} for similar recent works  based on  $(p,q)$-integers in the field of approximation theory. \\

In this paper, we consider  $0<q<p \leq 1$ and a King type modification of Lupa\c{s} Bernstein operators defined in \cite{lp} and investigate the statistical approximation properties of these operators. At last we show that this type of modification gives us better error estimation on some subintervals of $[0,1]$ than the classical  $(p,q)$-Lupa\c{s} Bernsrein operators. In case $p=1,$ it reduces to King type modification of $q$-Lupa\c{s} Bernstein operators. \\

\section{Construction of Operators}
 Now, we construct the King type modification of $(p,q)$-Lupa\c{s} Bernsrein operators (\ref{e9a}) which preserve monomials $e_i(x)= x^{i}$ for $(i= 0,2)$. For this study, we consider $0<q<p\leq 1$  satisfying the following condition\\
 \begin{equation}\label{bb}
 pq([n]-1) >p^{n}(p-q)
 \end{equation}
  for $n\geq 2.$\\

Let $r_{n}(x)$ be a sequence of real valued continuous functions defined on $[0,1]$ with $0\leq r_{n}(x)\leq1$. Let us consider the following operators:\\
		
\begin{equation}\label{18e}
L^*_{n}(f;p;q;x)= \sum\limits_{k=0}^{n}~~ \frac{f \bigg(\frac{p^{n-k}~[k]_{p,q}}{[n]_{p,q}}\bigg) ~\left[\begin{array}{c}
	n \\
	k%
	\end{array}%
	\right] _{p,q} p^{\frac{(n-k)(n-k-1)}{2}} q^{\frac{k(k-1)}{2}}~r_{n}(x)^{k}~(1-r_{n}(x))^{n-k}}{\prod\limits_{j=0}^{n-1}\{p^{j}(1-r_{n}(x))+q^{j} r_{n}(x)},
\end{equation}		
		
Where $f \in C[0,1] $, $ x \in [0,1]$ and $n\in N \backslash {0,1}$. it is clear that the operator $L^{*}_{n}(f,p,q;x)$ are positive and linear. Observe that, if we choose $r_{n}(x)=x$ then it turn out to be $(p,q)$-Lupa\c{s} Bernstein operators.

\begin{lemma}\label{b2}
	$L^*_{n}(f;p;q;x)$ satisfy the following properties.\\
	(1.)$L^*_{n}(e_0;p;q;x) = 1$\\
	(2.)$L^*_{n}(e_1;p;q;x) = r_{n}(x)$\\
	(3.)$L^*_{n}(f e_{2};p;q;x) = r_{n}^2(x) +\frac{x(1-x)p^{n-1}}{[n]_{p,q}}~-~\frac{r_{n}(x)^2(p-q)(1-r_{n}(x)}{p(1-r_{n}(x)+qr_{n}(x)}\big(1~-~\frac{p^{n-1}}{[n]_{p,q}}).$\\
\end{lemma}
\textbf{Note:} For our convenience, we denote $[n]_{p,q}~=~[n].$\\
		
Under the condition ( \ref{bb}), if we take\\

\begin{equation}\label{b3}
r_{n}(x)~=~-\frac{p^n+x^2[n](p-q)}{2(p^{n-1}q-p^{n}+q^2[n-1])}~+~\frac{\sqrt{{p^{2n}+x^{4}[n]^{2}(p-q)^{2}+2x^{2}[n](2pq([n]-1)-p^n(p-q))}}}{2(p^{n-1}q-p^{n}+q^2[n-1])},
\end{equation}		

then $L^*_{n}(f;p;q;x)$	preserve monomials, $L^*_{n}(e_0;p;q;x) = e_0(x) = 1$ and $L^*_{n}(e_2;p;q;x) = e_2 = x^2,$ for $n \geq 2$ .\\
Also $0 \leq r_{n}(x)\leq 1$ for $r_{n}(x)$ defined in (\ref{b3}).\\

From (\ref{bb}), we have\\
 \begin{equation}\label{cc}
2pq([n]-1)~-~p^{n}(p-q) > p^{n}(p-q)
 \end{equation}

 Using the inequality (\ref{cc}) we get\\

  \begin{equation}\label{c1}
p^{2n}~+~2x^{2}[n]p^{n}(p-q)~+~x^{4}[n]^{2}(p-q)^{2} = \big(p^{n}~+~x^{2}[n](p-q)\big)^{2}.
  \end{equation}

From above equality, we get $r_{n}(x)\geq 0$
 under the condition (\ref{bb}), since $(1-x)^2 \geq 0$ for $0 \leq x \leq 1,$ we have\\

   \begin{equation}\label{g1}
   p^{2n}~+~2x^{2}[n](2pq([n]-1)-p^{n}(p-q))~+~x^{4}[n]^{2}(p-q)^{2} \leq  \big(2(p^{n-1}q-p^{n}+q^2[n-1])+p^{n}~+~x^{2}[n](p-q)\big)^{2}.
   \end{equation}
    If we use (\ref{g1}) in (\ref{b3}) then we get $r_n(x)\leq 1$.\\
 		
\begin{rem}\label{r5.1}
	For $q\in(0,1)$ and $p\in(q,1]$, it is obvious that
	$\lim\limits_{n\to\infty}[n]_{p,q}=0 $ or $\frac1{p-q}$. In order to reach to convergence
	results of the operator $L^*_{n}(f;p;q;x),$ we take a sequence $q_n\in(0,1)$ and $p_n\in(q_n,1]$
	such that $\lim\limits_{n\to\infty}p_n=a,$ $\lim\limits_{n\to\infty}q_n=1$ and $\lim\limits_{n\to\infty}p_n^n=1,$ $\lim\limits_{n\to\infty}q_n^n=1$. So we get
	$\lim\limits_{n\to\infty}[n]_{p_n,q_n}=\infty$.
\end{rem}

\begin{thm}\label{1h1}
	Let $L^*_{n}(f;p;q;x)$ be the sequence of operators and the sequence $ p=p_n$ and $q=q_n$ satisfying Remark $\ref{r5.1}$  then for any function $f\in C[0,1]$\\
	$$\lim\limits_{n\to\infty}\lvert L^*_{n}(f;p;q;x_0)-f(x_0)\rvert =0$$\\
	for fixed $x_0\in [0,1]$\\
\end{thm}

\section{The Rates of Convergence}
		
The modulus of continuity for the space of function $ f\in C[0,1]$  is defined by\\

$$ w(f;\delta)=\sup\limits_{x,t\in C[0,1],~~ |t-x|<\delta} |f(t)-f(x)|$$\\

where ${w}(f;\delta)$ satisfies the following conditions:~~for all $f\in C[0,1],$\\

\begin{equation}\label{e118}
\lim\limits_{\delta\rightarrow0}~w (f;\delta) = 0.
\end{equation}
and
\begin{equation}\label{e119}
|f(t)-f(x)|\leq w(f;\delta)\bigg(\frac{|t-x|}{\delta}+ 1\bigg)		
\end{equation}				
		
Recall that, In \cite{iran} we obtained the following rate of convergence for the operators $(\ref{e1.1})$ for every $ f\in C[0,1]$ and $\delta> 0.$\\

\begin{equation}\label{e110}
|L_{n}(f;p;q;x)-f(x)|\leq w(f;\delta)\bigg(\frac{1}{\delta} \sqrt{{\frac{x(1-x)}{[n]}}}+ 1\bigg)		
\end{equation}

Now, we compute the rates of convergence of the operators $L^*_{n}(f;p;q;x)$ given by $(\ref{18e})$ to $f(x)$  by means of the modulus of continuity. and we also show that our error estimation is better than the $(p,q)$-Lupa\c{s} operator given by (\ref{e9a}).\\

\begin{theorem}\label{tt}
	Let $(p_n)$ and $(q_n)$ are the sequences satisfying remark (\ref{r5.1}) for each $n \geq 2$. For fixed $x \in [0,1],$ $f\in C[0,1]$ and $\delta_{n} >0$, we have\\
	
	$$|L^*_{n}(f,p_n,q_n;x)-f(x)|\leq 2w(f;\delta_{n}(x))$$\\
	Where\\
	
\begin{equation}
\delta_{n}(x)~=~  \sqrt {2x^{2} ~+~x\bigg(\frac{p^{n}_n+x^2[n](p_n-q_n)}{2(p^{n-1}_nq_n-p^{n}_n+q^{2}_n[n-1])}~-~\frac{\sqrt{{p^{2n}_n+x^{4}[n]^{2}(p_n-q_n)^{2}+2x^{2}[n](2p_nq_n([n]-1)-p^{n}_n(p_n-q_n))}}}{2(p^{n-1}_nq_n-p^{n}_n+q^{2}_n[n-1])}\bigg)}\tiny\\
\end{equation}
		
\end{theorem}

\begin{theorem}\label{2h2}
For all $f \in Lip_M(\rho)$\\

$\|L_{n}(f,p_n,q_n;x)-f(x)\|_{C[0,1]}~\leq~M\delta^{\rho}_n(x)$\\
		
where\\

$\delta_n(x) = \sqrt{\frac{x(1-x)}{[n]}}$
		
and M is a positive constant.
\end{theorem}

\begin{theorem}\label{3h3}
	For all $f \in Lip_M(\rho)$, under the condition (\ref{bb}) for $p~=~p_n$ and $q~=~q_n$, we have\\
	
$\|L^*_{n}(f,p_n,q_n;x)-f(x)\|_{C[0,1]}~\leq~M\delta^{\rho}_n(x)$\\
		
where\\
		
$\delta_n(x) =\sqrt{ 2x^2~+~x\bigg(\frac{p^n_n+x^2[n](p_n-q_n)}{2(p^{n-1}_nq_n-p^{n}_n+q^2_n[n-1])}~-~\frac{\surd{p^{2n}_n+x^{4}[n]^{2}(p_n-q_n)^{2}+2x^{2}[n](p^n(p_n-q_n))}}{2(p^{n-1}_nq_n-p^{n}_n+q^2_n[n-1])}\bigg)}$
		
and M is a positive constant.
\end{theorem}

\section{The Rates of Statistical Convergence}

	At this point, let us recall the concept of statistical convergence.
	The statistical convergence which was introduced by Fast \cite{LS1} in 1951, is an important research area in approximation theory. In \cite{LS1}, Gadjiev and Orhan used the concept of statistical convergence in approximation theory. They proved a Bohman-Korovkin type theorem for statistical convergence.\\
	Recently, statistical approximation properties of many operators are investigated \cite{bm,mab,mfa,bsm}.\\

	A sequence $x = (x_k)$ is said to be statistically convergent to a number L if for every $\epsilon > 0,$\\
	
	$$\delta \{K\in\textbf{ N}:|x_k -L|\geq \varepsilon\}~=~0,$$\\
	
	where $\delta(K)$ is the natural density of the set $K \subseteq \textbf{ N}.$\\
	
	The density of subset $K\subseteq N$ is defined by\\
	
	$$\delta(K):= \lim\limits_n\frac{1}{n}\{\text{the  number}~~ k \leq n : k\in K\}$$\\
	
	whenever the limit exists.\\

	For instance, $\delta(\textbf{N})~=~1, ~~\delta\{2K:k\in\textbf{ N}\}~=~\frac{1}{2}~~ \text{and}  ~~\delta\{k^2:K\in\textbf{ N}\}~=~0.$ \\
	
	To emphasize the importance of the statistical convergence, we have an example: The sequence\\
	\begin{equation}
	X_k =
	\left\{
	\begin{array}{lll}
	{L_1;~~~~~ if~~k = m^2}, & \\
	&\\
	{L_2;~~~~~~if~~k \neq m^2}. &  \\
	&\\
	\end{array}%
	\right.
	where~~ m\in\textbf{ N}
	\end{equation}
	
	is statistically convergent to $L_2$ but not convergent in ordinary sense when $L_1 \neq L_2.$ We note that any convergent sequence is statiscally convergent but not conversely. \\
	
	Now we consider sequences $ q = q_n$ and $p =p_n$ such that:\\
	
	\begin{equation}\label{1e1}
	st-\lim\limits_nq_n = 1,~~~~~~~st-\lim\limits_np_{n} = 1,~~~~~~~~ st-\lim\limits_nq^{n}_n = 1~~~~~~ \text{and}~~~~~~ st-\lim\limits_np^{n}_{n} = 1.
	\end{equation}
	
\begin{theorem}\label{4h4}
	If the sequences $p = p_n$ and $q= q_n$ satisfies the condition given in (\ref{1e1}), then\\
\begin{equation}\label{22}
|L^{*}_{n}(f,p_n.q_n;x)-f(x)|\leq 2w(f;\surd \delta_{n,x})\\
\end{equation}

for all $f \in C[0,1]$, where\\

$\delta_{n,x}~=~ L^{*}_{n}((t-x),^2p_n.q_n;x)$
\end{theorem}

\begin{theorem}\label{5h5}
If the sequences $p = p_n$ and $q= q_n$ satisfies the condition given in (\ref{1e1}), if $f\in C[0,1]$ then\\
	
$\|L_{n}(f,p_n,q_n;x)-f(x)\|_{C[0,1]}~\leq~2w(f;\delta_n)$\\

where\\
\begin{equation}\label{p6}
\delta_n~=~ \sqrt{\frac{2}{9[n]}}\\
\end{equation}

\end{theorem}

\begin{theorem}\label{6h6}
If the sequences $p = p_n$ and $q= q_n$ satisfies the condition given in (\ref{1e1}), if $f\in Lip_M(\rho)$ then\\

\begin{equation}\label{r2}
|L^{*}_{n}(f,p_n,q_n;x)-f(x)|~\leq~M\delta^{\rho}_n(x)\\
\end{equation}

where\\

$\delta_n(x) = \sqrt{\frac{2}{9[n]}}$\\

\end{theorem}


\end{document}